\documentclass[10pt]{article}
\usepackage{graphicx}
\usepackage{amsthm}
\usepackage{amsmath}
\usepackage{amssymb}

\begin{document}

\newtheorem{theorem}{Theorem}
\newtheorem{question}{Question}
\newtheorem{problem}{Problem}

\title{The Number of Halving Circles}
\author{Federico Ardila}
\date{}
\maketitle

\noindent {\bf 1. INTRODUCTION.} We say that a set $S$ of $2n+1$
points in the plane is in \emph{general position} if no three of
the points are collinear and no four are concyclic. We call a
circle \emph{halving} with respect to $S$ if it has three points
of $S$ on its circumference, $n-1$ points in its interior, and
$n-1$ in its exterior. The goal of this paper is to prove the
following surprising fact: \emph{any} set of $2n+1$ points in
general position in the plane has \emph{exactly} $n^2$ halving
circles.

\medskip

Our starting point is the following problem, which appeared in the
1962 Chinese Mathematical Olympiad \cite{Swetz}.

\begin{problem}
Prove that any set of $2n+1$ points in general position in the
plane has a halving circle.
\end{problem}

\noindent For the rest of sections 1 and 2, $n$ is a fixed
positive integer and $S$ signifies an arbitrary set of $2n+1$
points in general position in the plane.

There are several solutions to Problem 1. One possible approach is
the following. Let $A$ and $B$ be two consecutive vertices of the
convex hull of $S$. We claim that some circle going through $A$
and $B$ is halving. All circles through $A$ and $B$ have their
centers on the perpendicular bisector $\ell$ of the segment $AB$.
Pick a point $O$ on $\ell$ that lies on the same side of $AB$ as
$S$ and is sufficiently far away from $AB$ that the circle
$\Gamma$ with center $O$ and passing through $A$ and $B$
completely contains $S$. This can clearly be done. Now slowly
``push" $O$ along $\ell$, moving it towards $AB$. The circle
$\Gamma$ changes continuously with $O$. As we do this, $\Gamma$
stops containing some points of $S$. In fact, it loses the points
of $S$ one at a time: if it lost $P$ and $Q$ simultaneously, then
points $P,Q,A,$ and $B$ would be concyclic. We can move $O$
sufficiently far away past $AB$ that, in the end, the circle does
not contain any points of $S$.

Originally, $\Gamma$ contained all the points of $S$. Now, as it
loses one point of $S$ at a time in this process, we can decide
how many points we want it to contain. In particular, if we stop
moving $O$ when the circle is about to lose the $n$th point $P$ of
$S$, then the resulting $\Gamma$ is halving: it has $A$, $B$, and
$P$ on its circumference, $n-1$ points inside it, and $n-1$
outside it, as illustrated in Figure 1. $\Box$

\begin{figure}
\centering
\includegraphics[width=2in]{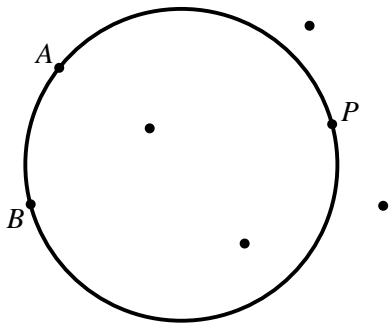}
\caption{A halving circle through $A, B,$ and $P$.}
\end{figure}

The foregoing proof shows that any set $S$ has several different
halving circles. We can certainly construct one for each pair of
consecutive vertices of the convex hull of $S$. In fact, the
argument can be modified to show that, for \emph{any} two points
of $S$, we can find a halving circle passing through them.

This suggests that we ask the following question: What can we say
about the number $N_S$ of halving circles of $S$? At first sight,
it seems that we really cannot say very much at all about this
number. Halving circles seem hard to ``control," and harder to
count.

We should, however, be able to find upper and lower bounds for
$N_S$ in terms of $n$. From the start we know that $N_S \geq
n(2n+1)/3$, since we can find a halving circle for each pair of
points of $S$, and each such circle is counted by three different
pairs. Computing an upper bound seems more difficult. If we fix
points $A$ and $B$ of $S$, it is indeed possible that all $2n-1$
circles through $A$, $B,$ and some other point of $S$ are halving.
The reader is invited to check this. Such a situation is not
likely to arise very often for a set $S$. However, it is not clear
how to make this idea precise, and then use it to obtain a
nontrivial upper bound.

For $n=2$, it is not too difficult to check by hand that $N_S=4$
for \emph{any} set $S$ of five points in general position in the
plane. This result first appeared in \cite{Perz}. It was also
proposed, but not chosen, as a problem for the 1999 International
Mathematical Olympiad. Notice that our lower bound gives $N_S \geq
4$.

In a different direction, a problem of the 1998 Asian-Pacific
Mathematical Olympiad, proposed by the author, asserted the
following.

\begin{problem}
$N_S$ has the same parity as $n$.
\end{problem}

Problem 2 follows easily from the nontrivial observation that, for
any $A$ and $B$ in $S$, the number of halving circles that go
through $A$ and $B$ is odd. We leave the proof of this observation
as a nice exercise.

\medskip

Amazingly, it turns out that we can say something much stronger.
The following result supercedes the previous considerations.

\begin{theorem}
Any set of $2n+1$ points in general position in the plane has
exactly $n^2$ halving circles.
\end{theorem}

Theorem 1 is the main result of this paper. In section 2 we prove
that every set of $2n+1$ points in general position in the plane
has the same number of halving circles. In section 3 we prove that
this number is exactly $n^2$, and we present a generalization.

\vspace{.5cm}

\noindent {\bf 2. THE NUMBER OF HALVING CIRCLES IS CONSTANT.} At
this point, we could cut to the chase and prove the very
counterintuitive Theorem 1. At the risk of making the argument
seem slightly longer, we believe that it is worthwhile to present
the motivation behind its discovery. Therefore, we ask the reader
to forget momentarily the punchline of this article.

Suppose that we are trying to find out whatever we can about the
number $N_S$ of halving circles of $S$. As mentioned in the
introduction, this number does not seem very tractable and it is
not clear how much we can say about it. Being optimistic, we might
hope to be able to answer the following two questions.

\begin{question}
What are the sharp lower and upper
bounds $m = m_{2n+1}$ and $M = M_{2n+1}$ for $N_S$?
\end{question}

\begin{question}
What are all the values that $N_S$
takes in the interval $[m, M]$?
\end{question}

Question 1 would appear to present considerable difficulty. To
answer it completely, we would first need to prove an inequality
$m \leq N_S \leq M$, and then construct suitable sets $S_{min}$
and $S_{max}$ that achieve these bounds. Let us focus on Question
2 instead. Here is a first approach.

Suppose that we start with the set $S_{min}$ (with $N_S=m$) and
move its points continuously so as to end up with $S_{max}$ (with
$N_S=M$). We might guess that the value of $N_S$ should change
``continuously," in the sense that $N_S$ should sweep out all the
integers between $m$ and $M$ as $S$ moves from a minimal to a
maximal configuration.

We know immediately that this would be overly optimistic. From
Problem 2 we learn that the parity of $N_S$ is determined by $n$,
so $N_S$ does not assume \emph{all} integral values between $m$
and $M$. In any case, the natural question to ask is: What kind of
changes does the value of $N_S$ undergo as $S$ changes
continuously?

\medskip

Let $S_{min}=\{P_1, \ldots, P_{2n+1}\}$ and $S_{max}=\{Q_1,
\ldots, Q_{2n+1}\}$. Now slowly transform $S_{min}$ into
$S_{max}$: first send $P_1$ to $Q_1$ continously along some path,
then send $P_2$ to $Q_2$ continuously along some other path, and
so on. We can think of our set $S$ as changing with time. At the
initial time $t=0$, our set is $S(0)=S_{min}$. At the final time
$t=T$, our set is $S(T)=S_{max}$. In between, $S(t)$ varies
continously with respect to $t$. Must $N_{S(t + \Delta t)} -
N_{S(t)}$ be small when $\Delta t$ is small? (As we move from
$S(0)$ to $S(T)$ continuously, it is likely that several
intermediate sets $S(t)$, with $0 < t < T$, are not in general
position. We shall see that we can go from $S(0)$ to $S(T)$ in
such a manner that we encounter only finitely many such sets. When
$S(t)$ is not in general position, we still need to know whether
$N_{S(t + \Delta t)} - N_{S(t - \Delta t)}$ must be small when
$\Delta t$ is small.)

In the way we defined the deformation from $S_{min}$ to $S_{max}$,
the points of $S$ move one at a time. Let us focus for the moment
on the interval of time during which $P_1$ moves towards $Q_1$.

Suppose that the number $N_S$ changes between time $t$ and time $t
+ \Delta t$. Then it must be the case that for some $i,j,k,$ and
$l$ the circle $P_iP_jP_k$ surrounds (or does not surround) point
$P_l$ at time $t$, but at time $t + \Delta t$ it does not (or
does) encircle $P_l$. For this to be true, it must happen that,
sometime between $t$ and $t + \Delta t$, either these four points
are concyclic or three of them are collinear. Since $P_1$ is the
only point that moves in this process, we can conclude that $P_1$
must cross a circle or a line determined by the other points; this
is what causes $N_S$ to change. We will call the circles and lines
determined by the points $P_2, P_3, \ldots, P_{2n+1}$ the
\emph{boundaries}.

We are free to choose the path along which $P_1$ moves towards
$Q_1$. To make things easier, we may assume that $P_1$ never
crosses two of the boundaries at the same time. This can clearly
be guaranteed: we know that these boundaries intersect pairwise in
finitely many points, and all we have to do is avoid their
intersection points in the path from $P_1$ to $Q_1$. We can also
assume that $\Delta t$ is small enough that $P_1$ crosses exactly
one boundary between times $t$ and $t+\Delta t$. Let us see how
$N_S$ can change in this time interval.

It will be convenient to call a circle $P_iP_jP_k$
\emph{$(a,b)$-splitting} (where $a+b = 2n-2$) if it has $a$ points
of $S$ inside it and the remaining $b$ points outside it. The
halving circles are the $(n-1,n-1)$-splitting circles.

\begin{figure}
\centering
\includegraphics[width=2.5in]{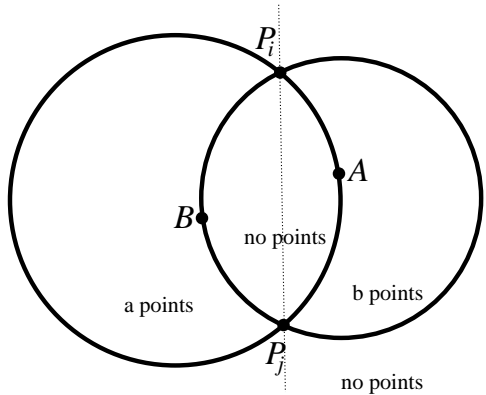}
\caption{$P_1$ crosses line $P_iP_j$.}
\end{figure}

Assume first that $P_1$ crosses line $P_iP_j$ in going from
position $P_1(t) = A$ to position $P_1(t+\Delta t) = B$. From the
remarks made earlier, we know that only circle $P_1P_iP_j$ can
change the value of $N_S$ by becoming or ceasing to be halving.
Assume that circle $AP_iP_j$ is $(a,b)$-splitting. Since $P_1$
only crosses the boundary $P_iP_j$ when going from $A$ to $B$, the
region common to circles $AP_iP_j$ and $BP_iP_j$ cannot contain
any points of $S$, as indicated in Figure 2. The region outside of
both circles cannot contain points of $S$ either. For circle
$AP_iP_j$ to be $(a,b)$-splitting, the other two regions must then
contain $a$ and $b$ points, respectively, as shown. Therefore
circle $BP_iP_j$ is $(b,a)$-splitting. It follows that $AP_iP_j$
is halving if and only if $BP_iP_j$ is halving (if and only if
$a=b=n-1$). Somewhat surprisingly, we conclude that the value of
$N_S$ does not change when $P_1$ crosses a line determined by the
other points; it can only change when $P_1$ crosses a circle.

Now assume that $P_1$ crosses circle $P_iP_jP_k$ in moving from
position $P_1(t) = A$ inside the circle to position $P_1(t+\Delta
t) = B$ outside it, as shown in Figure 3. (The other case, when
$P_1$ moves into the circle, is analogous.) The value of $N_S$ can
change only by circles $P_iP_jP_k$, $P_1P_jP_k$, $P_1P_kP_i$, and
$P_1P_iP_j$ becoming or ceasing to be halving. We can assume that
$P_1$ crosses the arc $P_iP_j$ of the circle that does not contain
point $P_k$. Notice that $A$ must be outside triangle $P_iP_jP_k$
if we want $P_1$ to cross only one boundary in the time interval
considered. Assume that circle $P_iP_jP_k$ is $(a,b)$-splitting
when $P_1 = A$. As before, we know that the only regions of Figure
3 containing points of $S$ are the one common to circles $AP_iP_j$
and $BP_iP_j$ and the one outside both of them. They must contain
$a-1$ and $b$ points, respectively. Circle $P_iP_jP_k$ goes from
being $(a,b)$-splitting to being $(a-1,b+1)$-splitting. The same
is true of circle $P_1P_iP_j$.

\begin{figure}
\centering
\includegraphics[width=3in]{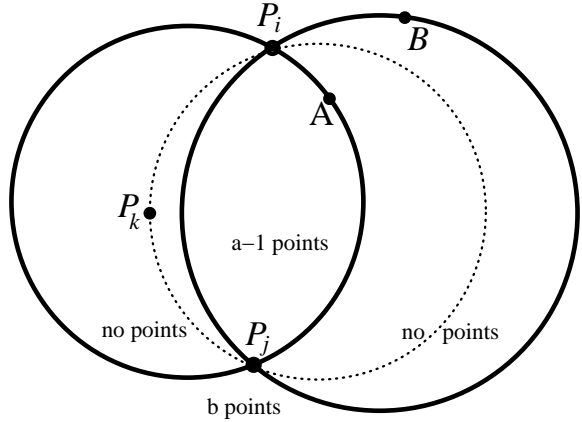}
\caption{$P_1$ crosses circle $P_iP_jP_k$.}
\end{figure}

It is also not hard to see, by a similar argument, that circles
$P_1P_jP_k$ and $P_1P_kP_i$ both change from being
$(a-1,b+1)$-splitting to being $(a,b)$-splitting. Again, the key
assumption is that $P_1$ only crosses the boundary $P_iP_jP_k$ in
this time interval.

So, by having $P_1$ cross circle $P_iP_jP_k$, we have traded two
$(a,b)$-splitting and two $(a-1,b+1)$-splitting circles for two
$(a-1,b+1)$-splitting and two $(a,b)$-splitting circles,
respectively. It follows that the number $N_S$ of halving circles
also remains constant when $P_1$ crosses a circle $P_iP_jP_k$!

We had shown that, as we moved $P_1$ to $Q_1$, $N_S$ could only
possibly change in a time interval when $P_1$ crossed a boundary
determined by the other points. But now we see that, even in such
a time interval, $N_S$ does not change! Therefore moving $P_1$ to
$Q_1$ does not change the value of $N_S$. Similarly, moving $P_i$
to $Q_i$ does not change $N_S$ for $i = 2, 3, \ldots, 2n+1$. It
follows that $N_S$ is the same for $S_{min}$ and $S_{max}$. In
fact, $N_S$ is the same for any set $S$ of $2n+1$ points in
general position!

\vspace{.5cm}

\noindent {\bf 3. THE NUMBER OF HALVING CIRCLES IS $n^2$.} Now
that we know that the number $N_S$ depends only on the number of
points in $S$, let $N_{2n+1}$ be the number of halving circles for
a set of $2n+1$ points in general position. We compute $N_{2n+1}$
recursively.

Construct a set $S$ of $2n+1$ points as follows. First consider
the vertices of a regular $(2n-1)$-gon with center $O$. Now move
them very slightly so that they are in general position. Label
them $P_1, \ldots, P_{2n-1}$ clockwise. The deformation should be
sufficiently slight that all the lines $OP_i$ still split the
remaining points into two sets of equal size, and all the circles
$P_iP_jP_k$ still contain $O$. Also consider a point $Q$ located
sufficiently far away from the others that it lies outside all the
circles formed by the points considered so far. Of course, we need
$Q$ to be in general position with respect to the remaining
points. We count the number of halving circles of $S=\{O, P_1,
\ldots, P_{2n-1}, Q\}$.

First consider the circles of the form $P_iP_jP_k$. These circles
contain $O$ and do not contain $Q$, so they are halving for $S$ if
and only if they are halving for $\{P_1, \ldots, P_{2n-1}\}$. Thus
there are $N_{2n-1}$ such circles.

Next consider the circles $OP_iP_j$. It is clear that these
circles contain at most $n-2$ other $P_k$s. They do not contain
$Q$, so they contain at most $n-2$ points, and they are not
halving.

Finally consider the circles that go through $Q$ and two other
points $X$ and $Y$ of $S$. Circle $QXY$ splits the remaining
points in the same way that line $XY$ does. More precisely, circle
$QXY$ contains a point $P$ of $S$ if and only if $P$ is on the
same side of line $XY$ that $Q$ is. This follows easily from the
fact that $Q$ lies outside circle $PXY$. Therefore we have to
determine which lines determined by two points of $S-\{Q\}$ split
the remaining points of this set into two subsets of $n-1$ points
each. This question is much easier to answer: the lines $OP_i$ do
this and the lines $P_iP_j$ do not. It follows that the $2n-1$
circles $OP_iQ$ are halving, and the circles $P_iP_jQ$ are not.

To summarize: the halving circles of $S$ are the $N_{2n-1}$
halving circles of $\{P_1, \ldots, P_{2n-1}\}$ and the $2n-1$
circles $OP_iQ$. Therefore $N_{2n+1}=N_{2n-1}+2n-1$. Since
$N_3=1$, it follows inductively that $N_{2n+1}=n^2.$ This
completes the proof of Theorem 1. $\Box$

\begin{theorem}
Consider a set of $2n+1$ points in general position in the plane,
and two nonnegative integers $a$ and $b$ satisfying $a<b$ and
$a+b=2n-2$. There are exactly $2(a+1)(b+1)$ circles that are
either $(a,b)$-splitting or $(b,a)$-splitting.
\end{theorem}

\noindent \emph{Sketch of proof.} The argument of section 2
carries over directly to this situation and shows that the number
of circles under consideration, which we denote $N(a,b)$, depends
only on $a$ and $b$. Therefore, it suffices to compute it for the
set $S$ constructed in the proof of Theorem 1.

Just as earlier, there are $N(a-1,b-1)$ such circles among the
circles $P_iP_jP_k$. Among the $OP_iP_j$ there are exactly $2n-1$
such circles, namely, the circles $OP_iP_{i+a+1}$ (taking
subscripts modulo $2n-1$). There are also $2n-1$ such circles
among the $QP_iP_j$, namely, the circles $QP_iP_{i+a+1}$. Finally,
there are no such circles among the $OP_iQ$. Therefore
\[
N(a,b) = N(a-1,b-1) + 4n - 2 = N(a-1,b-1) + 2a + 2b + 2.
\]
For $a=0$, we get that $N(0,b)=2b+2$. Theorem 2 then follows by
induction. $\Box$

\medskip

It is worth mentioning that our study is closely related to the
Voronoi diagram and the Delaunay triangulation of a point
configuration. The language of oriented matroids provides a very
nice explanation of this connection; for details, see \cite[sec.
1.8]{Bjorner}. In fact, Theorems 1 and 2 are essentially
equivalent to a beautiful result of D.T. Lee \cite{Lee}, which
gives a sharp bound for the number of vertices of an order $j$
Voronoi diagram. See \cite[Theorem 3.5]{Clarkson} for another
proof.

Under a stereographic projection, the halving circles of a point
configuration in the plane correspond to the halving planes of a
point configuration on a sphere in three-dimensional space. More
generally, we could also attempt to count the halving hyperplanes
of a point configuration in $n$-dimensional space. This problem
belongs to the vast literature on $k$-sets and $j$-facets, where
exact enumerative results are very rare. As an introduction, we
recommend \cite[chap. 11]{Matousek} to the interested reader.

\vspace{.5cm}

\noindent {\bf 4. ACKNOWLEDGEMENTS.} The results in this article
were first presented by the author at Dan Kleitman's birthday
conference at the Massachusetts Institute of Technology in 1999.
The author wishes to thank Dan Kleitman, Richard Stanley, Timothy
Chow and Uli Wagner for useful discussions on this subject during
and after the conference.

\medskip

\small

\noindent {\bf FEDERICO ARDILA} was born in Bogot\'{a}, Colombia,
in 1977. He received his B.Sc. and Ph.D. degrees in mathematics
from the Massachusetts Institute of Technology. He is currently a
postdoctoral researcher at the Mathematical Sciences Research
Institute. His research interests include geometric, algebraic,
and enumerative combinatorics. He is also closely involved with
the Colombian and International Mathematical Olympiads, through
which he first found an interest in mathematics; he currently
serves on the IMO Advisory Board. He likes to spend his free time
exploring Bogot\'{a} (when he is lucky enough to be there), on the
soccer field, or treasure hunting in little music stores.

\noindent \emph{Mathematical Sciences Research Institute,
Berkeley, CA 94720, USA \\} \emph{federico@msri.org}
\end{document}